\documentclass[12pt]{article}

\usepackage{latexsym,amssymb}
\usepackage{psfig}

\def\tr{^{\rm T}}
\def\Real{\mathbb R}

\def\qmx#1{\left(\matrix{#1}\right)}

\def\beeq#1{\begin{equation}{#1}\end{equation}}

\def\ba{\begin{array}}
\def\ea{\end{array}}
\def\eqa{\begin{eqnarray}}
\def\eqe{\end{eqnarray}}

\newtheorem{proposition}{Proposition}

\newtheorem{lemma}{Lemma}

\newenvironment{proof}{\medskip\noindent{\it Proof. }}{ \medskip}
\newenvironment{remark}{\medskip\noindent{\it Remark. }}{
\medskip}

\begin{document}

\title{Nonlinear Internal Models for Output Regulation \thanks{This work was partially
supported by NSF under grant ECS-0314004, by the AFOSR under grant
F49620-01-10039, and by the Boeing-McDonnell Douglas Foundation.}}

\author{C.I.Byrnes $^{\dag}$, A.Isidori $^{\dag \ddag}$}
\date{}

\maketitle
\begin{center}

$^{\dag}$Department of Systems Science and Mathematics, Washington
University, \\St. Louis, MO 63130.

\smallskip
$^{\ddag}$Dipartimento di Informatica e Sistemistica,
 Universit\`{a} di Roma ``La Sapienza'', \\00184 Rome, ITALY.
\end{center}
\maketitle

\vskip0.5in
\begin{abstract}
In this paper we show how nonlinear internal models can be
effectively used in the design of output regulators for nonlinear
systems. This result provides a significant enhancement of the
non-equilibrium theory for output regulation, which we have
presented in the recent paper ``Limit Sets, Zero Dynamics, and
Internal Models in the Problem of Nonlinear Output Regulation".
\end{abstract}

\noindent{\small {\bf Keywords}: Limit Sets, Zero Dynamics,
Internal Model, Regulation, Tracking, Nonlinear Control.}

\section{Introduction}

In the recent paper \cite{BI03}, we have laid the foundations for
a non-equilibrium theory of nonlinear output regulation, giving a
more general (non-equilibrium) definition of the problem, deriving
necessary conditions for the existence of solutions, and
describing how the necessary conditions thus found, complemented
with an additional set of appropriate hypotheses, can be used in
the design of output regulators.

We recall that, generally speaking, the problem of output
regulation is to have the regulated variables of a given
controlled plant to asymptotically track (or reject) all desired
trajectories (or disturbances) generated by some fixed autonomous
system, known as exosystem. The hypotheses assumed in \cite{BI03}
for the design of output regulators no longer include the
assumption, common to all earlier literature, that the
zero-dynamics of the controlled plant have a globally
asymptotically stable equilibrium. Rather, this assumption is
replaced with the (substantially weaker) hypothesis that the zero
dynamics of the plant ``augmented by exsosystem" have a compact
attractor. In \cite{BI03}, though, we have retained the (rather
strong) assumption, itself also common to all earlier literature,
that the set of all ``feedforward inputs capable to secure perfect
tracking" is a subset of the set of solutions of a suitable {\em
linear} differential equation. In this note we show that, within
the new framework, the assumption of linearity can also be
dropped.

Since this technical note is to be viewed as a continuation of the
work \cite{BI03}, we retain the same notation and -- to avoid
duplication -- we refer the reader to a number key concepts
introduced and/or described in that paper, among which the notion
of {\em omega limit set $\omega(B)$ of a set $B$}, plays a major
role. This concept is a deep generalization of the classical
concept, due to Birkhoff, of omega limit set of a point and
provides a rigorous definition of steady-state response in a
nonlinear system (see \cite{BI03} for details).

\section{Problem Statement}

In \cite{BI03}, as an illustration of how the new non-equilibrium
concepts can be applied to the design of regulators, we have shown
how the problem of output regulation can be solved, under
appropriate assumptions, for a system which can be put in the form
\beeq{\label{pla1}\ba{rcl} \dot z &=& f_0(z,w) +
f_1(z,\zeta,w)\zeta
\\ \dot \zeta &=& q(z,\zeta,w) + u\\ e &=&\zeta\\
y &=& \zeta\,,\ea} with state $(z,\zeta)\in \Real^n\times \Real$,
control input $u\in \Real$, regulated output $e\in \Real$,
measured output $w\in \Real$ and exogenous (disturbance) input
$w\in \Real^r$ generated by an exosystem \beeq{\label{exo} \dot w
= s(w)\,.} The functions $f_0(z,w),f_1(z,\zeta,w),q(z,\zeta,w)$
are $C^k$ functions (for some large $k$) of their arguments.

\begin{remark} System (\ref{pla1}) may look very particular, as it
has relative degree 1 between control input $u$ and regulated
output $e$. However, the design methodology described in
\cite{BI03}, and pursued in what follows under much weaker
hypotheses, lends itself to a straightforward extension to systems
with higher relative degree, namely systems having the form of
equation (33) in \cite{BI03}. Details are somewhat lengthy and for
this reason they are not included here. $\triangleleft$
\end{remark}

The analysis in \cite{BI03} was based on three standing
hypotheses. The first of these hypotheses is that the exosystem is
{\em Poisson stable}, namely that:

\medskip\noindent {\em Assumption} 0. The set $W \subset \Real^r$ of admissible initial conditions for the exosystem
(\ref{exo}) is {\em compact}  and $W=\bigcup_{w\in W}\omega(w)$.
$\triangleleft$

\medskip

Letting  $Z\times E \subset \Real^n \times \Real$ be the {\em
compact} set of initial states of (\ref{pla1}) for which the
problem of output regulation is to be solved, the second
hypothesis is that the trajectories of the zero dynamics of
(\ref{pla1}), augmented with (\ref{exo}), are {\em bounded},
namely that:

\medskip
\noindent {\em Assumption} 1. The positive orbit of $Z\times W$
under the flow of \beeq{\label{zerodyn} \ba{rcl} \dot z &=&
f_0(z,w) \\ \dot w &=&s(w)\ea} has a compact closure, and
$\omega(Z\times W) \subset {\rm int}(Z)\times W$. $\triangleleft$

If follows from this assumption that the set \[ {\cal A}_0 :=
\omega(Z\times W)\,,\] i.e the $\omega$-limit set  -- under the
flow of (\ref{zerodyn}) -- of the set $Z\times W$,  is a nonempty,
compact, invariant set which is stable in the sense of Lyapunov
and uniformly attracts $Z\times W$.

The third assumption was that:

\medskip
\noindent {\em Assumption} 2. There exist an integer $d$ and real
numbers $a_0, a_1, \ldots, a_{d-1}$ such that, for any
$(z_0,w_0)\in {\cal A}_0$, the solution $(z(t),w(t))$ of
(\ref{zerodyn}) passing through $(z_0,w_0)$ at $t=0$ is such that
the function $\varphi(t):=-q(z(t),0,w(t))$ satisfies
\[
\varphi^{(d)} + a_{d-1}\varphi^{(d-1)}+ \cdots +a_1\varphi^{(1)}+
a_0\varphi =0\;. \;\triangleleft
\]

It is well-known that the function $\varphi(t)$ considered above
is the input needed to keep the regulated output $e(t)$
identically zero (so long as $e(0)=0$). The assumption above calls
for the existence of a {\em linear} differential equation of which
$\varphi(t)$ be a solution. In this paper, we drastically weaken
this assumption, by simply calling for the existence of a {\em
nonlinear} differential equation of which $\varphi(t)$ be a
solution, namely:

\medskip
\noindent{\em Assumption} 2-nl. There exists an integer $d$ and a
locally Lipschitz function $f: \Real^d \to \Real$ such that, for
any $(z_0,w_0)\in {\cal A}_0$, the solution $(z(t),w(t))$ of
(\ref{zerodyn}) passing through $(z_0,w_0)$ at $t=0$ is such that
the function $\varphi(t):=-q(z(t),0,w(t))$ satisfies
\beeq{\label{immnonlin} \varphi^{(d)} + f(\varphi,\varphi^{(1)},
\ldots, \varphi^{(d-1)})=0\;. \;\triangleleft }

\section{Output Regulation via Nonlinear Internal Models}

We proceed now with the construction of a controller which solves
the problem of output regulation for system (\ref{pla1}). To this
end, consider the sequence of functions recursively defined as
\[\tau_1(z,w)=-q(z,0,w)\,,\quad \ldots , \quad \tau_{i+1}(z,w) =
{\partial \tau_i \over \partial z}f_0(z,w) + {\partial \tau_i
\over \partial w}s(w)\,,
\]
for $i=1, \ldots, d-1$, and consider the map
 \[ \ba{rcccl}
 \tau &:& \Real^n\times \Real^r &\to & \Real^d\\[1mm]

 &&(z,w) &\mapsto& {\rm
col}(\tau_1(z,w), \tau_2(z,w), \ldots, \tau_d(z,w))\,.\ea\] If
$k$, the degree of continuous differentiability of the functions
in (\ref{pla1}), is large enough, the map  $\tau$ is well defined
and $C^1$. In particular $\tau({\cal A}_0)$, the image of
 ${\cal A}_0$ under $\tau$\, is a {\em compact}
subset of $\Real^d$, because ${\cal A}_0$ is a compact subset of
$\Real^n\times \Real^r$.

Let $ f_{\rm c}: \Real^d \to \Real$ be any locally Lipschitz
function of compact support which agrees  on $\tau({\cal A}_0)$
with the function $ f$ defined in (\ref{immnonlin}), i.e. a
function such that, for some compact superset ${\cal S}$ of
$\tau({\cal A}_0)$ satisfies
\[\ba{rcll} f_{\rm c}(\eta) &=& 0 & \mbox{for all $\eta \in {\cal S}$} \\[1mm] f_{\rm c}(\eta) &=& f(\eta) & \mbox{for all $\eta \in
\tau({\cal A}_0)$}.\ea\] Note that, by definition, there is a
number $C>0$ such that \beeq{\label{bound} | f_{\rm c}(\eta)| \le
C\,, \qquad \mbox{for all $\eta \in \Real^d$}. }

Let now $\Phi_{\rm c}(\cdot)$ be the smooth vector field of
$\Real^d$ defined as
\[
\Phi_{\rm c}(\eta) = \qmx{\eta_2 \cr \eta_3 \cr \cdots  \cr \eta_d
\cr -f_{\rm c}(\eta_1,\eta_2,\ldots,\eta_d)\cr}
\]
and let $\Gamma$ be the $\Real^{1\times d}$ matrix
\[
\Gamma = \qmx{1 & 0 & \cdots & 0\cr}\,.
\]

Trivially, for any $(z_0,w_0)\in {\cal A}_0$ the function
$\tau(z(t),w(t))$ is such that
\[
{{\rm d} \over {\rm d}t}\,\tau_i(z(t),w(t))= \tau_{i+1}(z(t),w(t))
\]
for $i=1,\ldots, d-1$. If, in addition, Assumption 2-nl holds,
then
\[
{{\rm d} \over {\rm d}t}\,\tau_d(z(t),w(t))= -f(\tau(z(t),w(t))) =
-f_{\rm c}(\tau(z(t),w(t)))\,.
\]
Thus, if Assumption 2-nl holds, for any $(z_0,w_0)\in {\cal A}_0$
the function $\tau(z(t),w(t))$ satisfies \beeq{\label{imm1} {{\rm
d} \over {\rm d}t}\,\tau(z(t),w(t)) = \Phi_c(\tau(z(t),w(t)))\;.}

\begin{remark} In other words, the restriction of (\ref{zerodyn}) to ${\cal
A}_0$, which is invariant under the flow of (\ref{zerodyn}), and
the restriction of \beeq{\label{imm1eta}\dot \eta = \Phi_c(\eta)}
to $\tau({\cal A}_0)$, which is invariant under the flow of
(\ref{imm1eta}), are $\tau$-related systems. $\triangleleft$
\end{remark}

Moreover, \beeq{\label{imm2} q(z,0,w)=-\Gamma\tau(z,w)\,, \qquad
\mbox{for all $(z,w)\in {\cal A}_0.$}}

\bigskip
Consider now a control law of the form
\beeq{\label{contr1}\ba{rcl} \dot \xi &=& \Phi_{\rm c}(\xi) + Gv\\
u &=& \Gamma\xi + v\,, \ea}  with \beeq{\label{contr2} v = -
ky\;.} A simple check shows that, if Assumptions 0, 1 and 2-nl
hold, this controller does have the {\em internal-model property}
(see \cite{BI03}), relative to the set ${\cal A}_0$.

The control of (\ref{pla1}) by means of (\ref{contr1}) results in
a system \beeq{\label{closprovv}\ba{rcl} \dot z &=&
f_0(z,w)+f_1(z,\zeta,w)\zeta\\ \dot \zeta &=& q(z,\zeta,w) +
\Gamma\xi + v
\\ \dot \xi &=& \Phi_{\rm c}(\xi) + Gv\\ \dot w &=&s(w)\,,\ea}
which, viewing $v$ as input and $\zeta$ as output,  has relative
degree 1 and a zero-dynamics
\[\ba{rcl} \dot z &=&
f_0(z,w)
\\ \dot \xi &=& \Phi_{\rm c}(\xi)+G(-q(z,0,w)-\Gamma \xi)\\ \dot w &=&s(w)\,.\ea\]

The asymptotic properties of the latter are summarized in the
following results.

\begin{lemma}\label{LM6.2}
Suppose Assumptions {\rm 0, 1} and {\rm 2-nl} hold.  Consider the
triangular system \beeq{\ba{rcl} \label{aux1}\dot z &=&
f_0(z,w)\\\dot w &=& s(w)\\ \dot \xi &=& \Phi_{\rm
c}(\xi)+G(-q(z,0,w)-\Gamma \xi)\,.\ea} Let $\Xi \subset \Real ^d$
be a bounded subset such that $\tau({\cal A}_0)\subset {\rm
int}(\Xi)$. Choose numbers $c_0, c_1, \ldots, c_{d-1}$ such that
the polynomial \[ p(\lambda) = \lambda^d + c_{d-1}\lambda^{d-1} +
\cdots c_1\lambda +c_0\] is Hurwitz, set
\[\ba{rcl}
G_0&=& {\rm col}(c_0, c_1, \ldots,
c_{d-1})\\[2mm]
D_\kappa&=&{\rm diag}(\kappa, \kappa^2, \ldots ,\kappa^d)\ea
\]
and \[
 G = D_\kappa G_0\,.
 \] Then, there is a number $\kappa^\ast$ such that, if $\kappa\ge
\kappa^\ast$,
\[{\rm graph}(\tau\vert_{{\cal A}_0})=\omega(Z\times W\times
\Xi)\] under the flow of (\ref{aux1}). In particular, ${\rm
graph}(\tau\vert_{{\cal A}_0})$ is a compact invariant set which
uniformly attracts $Z\times W\times \Xi$.
\end{lemma}

\begin{proof}  Set,
for convenience, $B=Z\times W\times \Xi$.

\medskip  First of all, we  show that the positive orbit of $B$ is
bounded, so that $\omega(B)$ is nonempty, compact, invariant and
uniformly attracts $B$ (see Lemma 2.1 of \cite{BI03}). To this
end, set
\[
\bar \xi= D_{\kappa}^{-1}\xi\,,
\]
and observe that the variable thus defined obeys
\beeq{\label{dote} \dot {\bar \xi} = \kappa (A-G_0\Gamma){\bar
\xi} + D_{\kappa}^{-1}Bf_{\rm c}(D_\kappa{\bar \xi})- G_0q(z,0,w)}
in which
\[
A = \qmx{0&1&\cdots &0&0\cr0&0&\cdots &0&0\cr \cdot&\cdot&\cdots
&\cdot&\cdot\cr 0&0&\cdots &0&1\cr 0&0&\cdots &0&0\cr}, \qquad
B=\qmx{0\cr 0\cr \cdot \cr 0\cr 1\cr}.\]

It is seen from (\ref{bound}) that, for any $\kappa >1 $ (which we
can assume without loss of generality) and for any $\bar \xi \in
\Real^d$,
\[|D_{\kappa}^{-1}Bf_{\rm c}(D_\kappa{\bar \xi})|\le C\,.
\]
Thus, the integral curves of (\ref{dote}) are maximally defined on
the entire real axis. Knowing that $\xi(t)$ exists for any $t>0$,
and setting
 \[
u(t) = D_{\kappa}^{-1}Bf_{\rm c}(D_\kappa{\bar
\xi}(t))-G_0q(z(t),0,w(t))\,,
\]
system (\ref{dote}) can be rewritten in the form
\[
\dot {\bar \xi} = \kappa (A-G_0\Gamma){\bar \xi} +u\,.\] Since the
matrix $(A-G_0\Gamma)$ is Hurwitz by construction, and $u(t)$ is
bounded (by some fixed number which only depends on the choice of
$W$ and $ Z$),  it is concluded that also $|\bar \xi(t)|$ is
bounded. As a consequence, $\xi(t)$ itself is bounded, by number
which depends on the choice of $W, Z$ and $\Xi$ and of the
parameter $\kappa$. This concludes the proof that the positive
orbit of $B$ is bounded.

\medskip  To find $\omega(B)$, consider now the (closed) solid
cylinder ${\cal C} = {\cal A}_0 \times \Real^d$. A simple argument
shows that $\omega(B) \subset {\cal C}$. In fact, for any point
$p=(z,w,\xi) \in \omega(B)$, there must exist sequences
$p_k=(z_k,w_k,\xi_k)$ and $t_k$, with $t_k\to \infty$ as $k\to
\infty$, such that $ \lim_{k \to \infty}\phi(t_k,p_k) = p$.
Because of the triangular structure of (\ref{aux1}), the point
$(z,w)$ is by definition in ${\cal A}_0$ and hence $(z,w,\xi)\in
{\cal C}$.

\medskip
It is immediate to see that ${\rm graph}(\tau\vert_{{\cal A}_0})$
is invariant for (\ref{aux1}). In fact, observe that, for any
$(z_0,w_0) \in {\cal A}_0$, the solution $(z(t),w(t))$ of the top
two equations of (\ref{aux1}) passing through $(z_0,w_0)$ at time
$t=0$ remains in ${\cal A}_0$ for all $t\in \Real$, because the
latter is invariant for (\ref{zerodyn}). Thus, bearing in mind
(\ref{imm2}), for all such trajectories the bottom equation of
(\ref{aux1}) can be rewritten as
\[
\dot \xi = \Phi_{\rm c}(\xi) + G(\Gamma\tau(z,w)-\Gamma\xi)\,.
\]
 Set $\chi(t) = \xi(t) - \tau(z(t),w(t))$.
Then, using (\ref{imm1}), we have \beeq{\label{dotchi} \dot \chi =
\Phi_{\rm c}(\chi+\tau(z,w))- \Phi_{\rm
c}(\tau(z,w))-G\Gamma\chi\,. } The point $\chi=0$ is an
equilibrium of this equation and, therefore, ${\rm
graph}(\tau\vert_{{\cal A}_0})$ is invariant for (\ref{aux1}). In
particular, for any point $p\in {\rm graph}(\tau\vert_{{\cal
A}_0})$ it is possible to find sequences $p_k\in {\rm
graph}(\tau\vert_{{\cal A}_0})$ and $t_k$, with $t_k\to \infty$ as
$k\to \infty$, such that $ \lim_{k \to \infty}\phi(t_k,p_k) = p$.
This shows that ${\rm graph}(\tau\vert_{{\cal A}_0})\subset
\omega(B)$.

\medskip To complete the proof, we need to show that no point in
${\cal C}\setminus {\rm graph}(\tau\vert_{{\cal A}_0})$ can be a
point of $\omega(B)$. For, pick any point $p_0 \in {\cal
C}\setminus {\rm graph}(\tau\vert_{{\cal A}_0})$, i.e. a point
$p_0=(z_0,w_0,\xi_0)$ in which $(z_0,w_0)\in {\cal A}_0$ and
$\xi_0 \ne \tau(z_0,w_0)$. The complete orbit of (\ref{aux1})
through $p_0$ is in ${\cal C}$, because ${\cal A}_0$  is invariant
under (\ref{zerodyn}). Thus, the $\xi$-component of the
corresponding trajectory is such that
$\chi(t)=\xi(t)-\tau(z(t),w(t))$ obeys (\ref{dotchi}), with
$\chi(0)\ne 0$.

Set
\[
\bar \chi= D_{\kappa}^{-1}\chi\,,
\]
and observe that the variable thus defined obeys
\beeq{\label{dote2} \dot {\bar \chi} = \kappa (A-G_0\Gamma){\bar
\chi} + B\Delta({\bar \chi},\tau(z,w),\kappa)\,} with
\[ \Delta({\bar \chi},\tau(z,w),\kappa)) = {1\over
\kappa^d}\Bigl[f_{\rm c}(D_\kappa {\bar \chi} + \tau(z,w))-f_{\rm
c}(\tau(z,w))\Bigr].
\]

Since the function $f_{\rm c}(\cdot)$ is locally Lipschitz, and
bounded, there exists a number $L$ such that, for any $\xi_1,\xi_2
\in \Real^d$,\[ |f_{\rm c}(\xi_1+\xi_2)-f_{\rm c}(\xi_2)|\le
L|\xi_1|\,.\] Thus, for any $\kappa >1$,
\[
|\Delta({\bar \chi},\tau(z,w),\kappa))|\le {1\over
\kappa^{\,d}}\,L|D_\kappa {\bar \chi}|\le L|{\bar \chi}|\,.
\]

On the other hand, the matrix $(A-G_0\Gamma)$ is Hurwitz by
construction and hence there exists a positive definite matrix $P$
such that
\[
P(A-G_0\Gamma)+ (A-G_0\Gamma)\tr P = -I\,.
\]
From this, it is  immediately seen that the function $V(\bar \chi)
= \bar \chi\,\tr P \bar \chi$ satisfies, along the trajectories of
(\ref{dote2}),
\[ \dot V(\bar \chi) \le -(\kappa - 2L|P|)|\bar \chi|^2\,.
\]
If $\kappa > 2L|P|$, there exist numbers $M>0$ and $\alpha>0$,
both depending on $\kappa$, such that \beeq{\label{expdecay}
|\chi(t)|\le Me^{-\alpha t}| \chi(0)|\,.}

This property can be used to show that no point of ${\cal
C}\setminus {\rm graph}(\tau\vert_{{\cal A}_0})$ can be a point of
$\omega(B)$. For, suppose by contradiction that a point
$p_0=(z_0,w_0,\xi_0)$, with $(z_0,w_0)\in {\cal A}_0$ and $\xi_0
\ne \tau(z_0,w_0)$, be a point of $\omega(B)$. Set
\beeq{\label{dzero} d_0 = |\xi_0 - \tau(z_0,w_0)| \,.}
 As the
trajectory of (\ref{aux1}) through this point is defined for all
$t\in \Real$ and bounded, there exists a number $K$ such that
$|\chi(t)|\le K$ for all $t<0$. Pick now a number $T>0$ such that
\[
Me^{-\alpha T}K \le 0.5 \,d_0\,.
\]
Then, using (\ref{expdecay}), we obtain
\[
|\xi_0 - \tau(z_0,w_0)| \le Me^{-\alpha T}|\xi(-T) -
\tau(z(-T),w(-T))| \le Me^{-\alpha T}K \le 0.5 \,d_0\,, \] which
contradicts (\ref{dzero}). This concludes the proof of the Lemma.
$\triangleleft$.

\end{proof}

\begin{lemma}\label{LM6.3}
Suppose Assumptions {\rm 0, 1} and {\rm 2-nl} hold.  Choose $G$ as
indicated in Lemma \ref{LM6.2}, with $\kappa \ge \kappa^\ast$. If
${\cal A}_0$ is locally exponentially attractive for
(\ref{zerodyn}), then ${\rm graph}(\tau\vert_{{\cal A}_0})$ is
locally exponentially attractive for (\ref{aux1}).
\end{lemma}

\begin{proof} To say that ${\cal A}_0$ is locally exponentially
attractive for (\ref{zerodyn}) is to say that there exists
positive numbers $\rho, A, \lambda$ such that, for any $(z_0,w_0)$
in the (closed) set
\[
{\cal A}_\rho = \{(z,w)\in \Real^n\times W: {\rm dist}((z,w),{\cal
A}_0)\le \rho\}\] the solution $(z(t),w(t))$ of (\ref{zerodyn})
passing through $(z_0,w_0)$ at $t=0$ is such that
\[
{\rm dist}((z(t),w(t)),{\cal A}_0) \le Ae^{-\lambda t}{\rm
dist}((z_0,w_0),{\cal A}_0)\,.
\]
Set, as in the proof of the previous Lemma, $\chi(t) =
\xi(t)-\tau(z(t),w(t))$. Then, \beeq{\label{dotchinew} \dot \chi =
\Phi_{\rm c}(\chi+\tau(z,w))- \Phi_{\rm c}(\tau(z,w))-G\Gamma\chi+
Gp(z,w)\,,}in which  \[ p(z,w) = \Gamma \tau(z,w) +q(z,0,w)\;.
\]
From (\ref{imm2}), it is seen that the function $p(z,w)$ vanishes
on ${\cal A}_0$. Since the function in question is $C^1$ and
${\cal A}_\rho$ is compact, there exist a number $N>0$ such that
\[
|p(z,w)|\le N{\rm dist}((z,w),{\cal A}_0), \qquad \mbox{for all
$(z,w)\in {\cal A}_\rho$ ,}\] and therefore, if ${\rm
dist}((z_0,w_0),{\cal A}_0)\le \rho^\prime$, with $\rho^\prime =
\rho/A$, \beeq{\label{expdecay1} |p(z(t),w(t))|\le NAe^{-\lambda
t}{\rm dist}((z_0,w_0),{\cal A}_0)\;. }

Changing $\chi$ into $\bar \chi=D_\kappa^{-1}\chi$ in
(\ref{dotchinew}) yields (see proof of the previous Lemma)
\[
\dot {\bar \chi} = \kappa (A-G_0\Gamma){\bar \chi} + B\Delta({\bar
\chi},\tau(z,w),\kappa)+G_0p(z,w)\] from which, using again the
same Lyapunov function $V(\bar \chi)$ introduced above and bearing
in mind (\ref{expdecay1}), standard arguments show that, if
$\kappa
> 2L|P|$, an estimate of the form
\[
|\chi(t)| \le \hat Me^{-\hat \alpha t}(|\chi(0)|+{\rm
dist}((z_0,w_0),{\cal A}_0))
\]
holds, for every $\chi(0)$ and for every $(z_0,w_0)$ satisfying
${\rm dist}((z_0,w_0),{\cal A}_0)\le \rho^\prime$. From this, the
result follows. $\triangleleft$ \end{proof}

As shown in \cite{BI03}, the properties indicated in these two
Lemma (which in \cite{BI03} were proven to hold under the much
strunger Assumption 2) guarantee that, if the input $v$ of
(\ref{closprovv}) is chosen as in (\ref{contr2}), the resulting
closed-loop system has the desired asymptotic properties.

As a matter of fact, the following conclusion holds.

\begin{proposition} Consider system (\ref{clos1}). Let  $Z, W$ be fixed compact
sets of initial conditions and suppose Assumptions {\rm 0, 1} and
{\rm 2-nl} hold. Pick any compact set $\Xi$ as indicated in Lemma
\ref{LM6.2} and any closed  interval $E$ of $\Real$. Choose $G$ as
indicated in Lemma \ref{LM6.2}, with $\kappa$ large enough so that
the conclusions of Lemma \ref{LM6.2} hold. Then, for every
$\varepsilon>0$, there exists a number $k^\ast$ such that, if
$k\ge k^\ast$,  the positive orbit of $Z \times W \times \Xi
\times E$ is bounded and there exists $\bar t$ such that
$|e(t)|\le \varepsilon$ for all $t \ge \bar t$. If, in addition,
${\cal A}_0$ is locally exponentially attractive for
(\ref{zerodyn}), then $e(t)\to 0$ as $t\to \infty$.
\end{proposition}

The proof of this Proposition follows a standard paradigm.
Changing the variable $\xi$ into $\eta=\xi-Ge$ and setting $\bar k
=k-\Gamma G$ yields a system of the form \beeq{\label{clos1}
\ba{rcl} \dot z &=& f_0(z,w) + f_1(z,e,w)e\\ \dot w &=& s(w) \\
\dot \eta &=& \Phi_{\rm c}(\eta+Ge) -G\Gamma(\eta + Ge) -
Gq(z,e,w) \\ \dot e &=& q(z,e,w) + \Gamma\eta - \bar ke\,. \ea}
This can be viewed as interconnection of two subsystems, one with
state $(z,w,\eta)$ and input $e$, the other with state $e$ and
input $(z,w,\eta)$. Then, small-gain arguments can be invoked to
show that, if $\bar k$ is large enough, the results of the
Proposition hold. Details on the proof that the subsystems in
question have the appropriate input-to-state stability properties
can be found in \cite{BIP}.


\begin{thebibliography}{99}
\bibitem{BI03} C.I. Byrnes, A. Isidori, Limit sets, zero dynamics and internal models
in the problem of nonlinear output regulation, {\em IEEE Trans. on
Automatic Control}, {\bf AC-48}, pp. 1712-1723, (2003).

\bibitem{BIP} C.I. Byrnes, A. Isidori, L. Praly, On the Asymptotic Properties of a System Arising in
Non-equilibrium Theory of Output Regulation, Preprint of the
Mittag-Leffler Institute, Stocknolm, June 2003.



\end{thebibliography}
\end{document}